\documentclass[10pt,oneside]{article}
\input amssymb.sty
\usepackage[cp1250]{inputenc}
\usepackage[T1]{fontenc}

\begin{document}
\newcommand{\qed}{\rule{1.5mm}{1.5mm}}
\newcommand{\proof}{\textit{Proof. }}
\newcommand{\ccon}{\rightarrowtail}
\newtheorem{theorem}{Theorem}[section]
\newtheorem{lemma}[theorem]{Lemma}
\newtheorem{remark}[theorem]{Remark}
\newtheorem{example}[theorem]{Example}
\newtheorem{corollary}[theorem]{Corollary}
\newtheorem{proposition}[theorem]{Proposition}
\newtheorem{claim}[theorem]{Claim}
\begin{center}
{\LARGE\textbf{Local approximation of the solutions of algebraic equations}\vspace*{3mm}}\\
{\large Marcin Bilski\footnote{E-mail address: Marcin.Bilski@im.uj.edu.pl\vspace*{1mm}\\
Partially supported by the grant NN201 3352 33 of the Polish
Ministry of Science and Higher Education.
}}\vspace*{2mm}\\
{\small\textit{Institute of Mathematics, Jagiellonian University,
Reymonta 4, 30-059 Krak\'ow,
Poland}}\vspace*{3mm}\\
\end{center}
{\small\textbf{Abstract.}We give a constructive proof of the fact
that every holomorphic solution\linebreak $y=f(x)$ of a system
$Q(x,y)=0$ of polynomial (or Nash) equations can be approximated
in a fixed neighborhood of every $x_0\in dom(f)$ by a sequence of
Nash
solutions. An algorithm of such approximation is next presented.\vspace*{2mm}\vspace*{0mm}}

\section{Introduction}\label{introduction}
The aim of this note is to present a simple geometric proof of the
following approximation theorem. The proof allows to construct an
algorithm which can be used in numerical computations (see Section
\ref{examp}).
\begin{theorem}\label{main}Let $U$ be an open subset of
$\mathbf{C}^n$ and let $f:U\rightarrow\mathbf{C}^k$ be a
holomorphic mapping that satisfies a system of equations
$Q(x,f(x))=0$ for $x\in U,$ where
$Q:\mathbf{C}^n\times\mathbf{C}^k\rightarrow\mathbf{C}^q$ is a
polynomial mapping. Then for every $x_0\in U $ there are an open
neighborhood $U_0\subset U$ and a sequence
$\{f^{\nu}:U_0\rightarrow\mathbf{C}^k\}$ of Nash mappings
converging uniformly to $f|_{U_0}$ such that $Q(x,f^{\nu}(x))=0$
for every $x\in U_0$ and $\nu\in\mathbf{N}.$
\end{theorem}

Local approximation of the solutions of algebraic or analytic
equations was investigated by M.~Artin in \cite{Ar1}, \cite{Ar2}
and \cite{Ar3} and Theorem~\ref{main} can be derived from the
results of these papers. Since our goal is to obtain an effective
procedure of approximation, the present note treats the problem in
a bit different way: first a reduction to the case where $Q$ is a
single polynomial (which, given the original equations, can be
computed) is carried out. The reduced problem is next solved and
here the combination of some of the ideas of \cite{Ar1} and
\cite{vD} is applied (see Section \ref{mainIprof}). The latter
article, due to L. van den Dries, deals with the global version of
the theorem for $f$ depending on one variable.

Theorem \ref{main} can be easily strengthened by replacing the
polynomial mapping $Q$ by a Nash mapping defined in some
neighborhood of the graph of $f$ (cf. Section \ref{mainIprof}).
The method of approximation presented in this note is more
efficient than the one we gave in \cite{B6}, where the question of
the existence of mappings approximating solutions of algebraic
equations is also discussed.

Both in the real case (M. Coste, J. M. Ruiz, M. Shiota \cite{CRS})
and in the complex case (L. Lempert \cite{Lem}) more general
(global) versions of Theorem \ref{main} are known to be true.
These results rely on the solution to the deep and important
M.~Artin's conjecture for which the reader is referred to
\cite{An}, \cite{Og}, \cite{Po1}, \cite{Po2}, \cite{Sp}. Such an
approach enabled the authors to reach the goal in an elegant and
relatively short way, but makes the proofs difficult to be applied
in effective computations.

Our interest in Theorem \ref{main} is partially motivated by
applications in the theory of analytic sets. In particular, papers
\cite{B1}--\cite{B3} contain results on approximation of complex
analytic sets by complex Nash sets whose proofs can be divided
into two stages: (i) preparation, where only direct geometric
methods appear, (ii) switching Theorem \ref{main}. Thus the
techniques of the present article allow to obtain local versions
of these facts in a purely geometric way.

The paper is organized as follows. Sections \ref{mainIprof} and
\ref{examp} are devoted to Theorem~\ref{main} and the algorithm
respectively. Preliminary material concerning Nash mappings and
sets as well as analytic sets with proper projection is gathered
in Section \ref{preli} below.
\section{Preliminaries}\label{preli}
\subsection{Nash mappings and sets}\label{prelnash}
Let $\Omega$ be an open subset of $\mathbf{C}^n$ and let $f$ be a
holomorphic function on $\Omega.$ We say that $f$ is a Nash
function at $x_0\in\Omega$ if there exist an open neighborhood $U$
of $x_0$ and a polynomial
$P:\mathbf{C}^n\times\mathbf{C}\rightarrow\mathbf{C},$ $P\neq 0,$
such that $P(x,f(x))=0$ for $x\in U.$ A holomorphic function
defined on $\Omega$ is said to be a Nash function if it is a Nash
function at every point of $\Omega.$ A holomorphic mapping defined
on $\Omega$ with values in $\mathbf{C}^N$ is said to be a Nash
mapping if each of its components is a Nash function.

A subset $Y$ of an open set $\Omega\subset\mathbf{C}^n$ is said to
be a Nash subset of $\Omega$ if and only if for every
$y_0\in\Omega$ there exists a neighborhood $U$ of $y_0$ in
$\Omega$ and there exist Nash functions $f_1,\ldots,f_s$ on $U$
such that $$Y\cap U=\{x\in U: f_1(x)=\ldots=f_s(x)=0\}.$$ We will
use the following fact from \cite{Tw}, p. 239. Let
$\pi:\Omega\times\mathbf{C}^k\rightarrow\Omega$ denote the natural
projection.
\begin{theorem}
\label{projnash} Let $X$ be a Nash subset of
$\Omega\times\mathbf{C}^k$ such that $\pi|_{X}:X\rightarrow\Omega$
is a proper mapping. Then $\pi(X)$ is a Nash subset of $\Omega$
and $dim(X)=dim(\pi(X)).$
\end{theorem}
The fact from \cite{Tw} stated below explains the relation between
Nash and algebraic sets.
\begin{theorem}Let $X$ be a Nash subset of an open set
$\Omega\subset\mathbf{C}^n.$ Then every analytic irreducible
component of $X$ is an irreducible Nash subset of $\Omega.$
Moreover, if $X$ is irreducible then there exists an algebraic
subset $Y$ of $\mathbf{C}^n$ such that $X$ is an analytic
irreducible component of $Y\cap\Omega.$
\end{theorem}

\subsection{Analytic sets}\label{setswithprop}
Let $U, U'$ be domains in $\mathbf{C}^n,\mathbf{C}^k$ respectively
and let
$\pi:\mathbf{C}^n\times\mathbf{C}^k\rightarrow\mathbf{C}^n$ denote
the natural projection. For any purely $n$-dimensional analytic
subset $Y$ of $U\times U'$ such that $\pi|_Y:Y\rightarrow U$ is a
proper mapping by $\mathcal{S}(Y,\pi)$ we denote the set of
singular points of $\pi|_{Y}:$
$$\mathcal{S}(Y,\pi)=Sing(Y)\cup\{x\in Reg(Y):(\pi|_Y)'(x) \mbox{ is not an isomorphism}\}.$$
We often write $\mathcal{S}(Y)$ instead of $\mathcal{S}(Y,\pi)$
when it is clear which projection is taken into consideration.

It is well known that $\mathcal{S}(Y)$ is an analytic subset of
$U\times U',$ $dim(Y)<n$ (cp. \cite{Ch}, p. 50), hence by the
Remmert theorem $\pi(\mathcal{S}(Y))$ is also analytic. Moreover,
the following hold. The mapping $\pi|_{Y}$ is surjective
and open and there exists an integer $s=s(\pi|_{Y})$ such that:\vspace*{2mm}\\
(1) $\sharp(\pi|_{Y})^{-1}(\{a\})<s$ for
$a\in\pi(\mathcal{S}(Y)),$\\
(2) $\sharp(\pi|_{Y})^{-1}(\{a\})=s$ for $a\in U\setminus
\pi(\mathcal{S}(Y)),$\\
(3) for every $a\in U\setminus \pi(\mathcal{S}(Y))$ there exists a
neighborhood $W$ of $a$ and holomor-\linebreak \hspace*{5.4mm}phic
mappings $f_1,\ldots,f_s:W\rightarrow U'$ such that $f_i\cap
f_j=\emptyset$ for $i\neq j$ and\linebreak
\hspace*{5.4mm}$f_1\cup\ldots\cup f_s=$$(W\times U')\cap
Y.$\vspace*{2mm}

Let $E$ be a purely $n$-dimensional analytic subset of $U\times
U'$ with proper projection onto a domain $U\subset\mathbf{C}^n,$
where $U'$ is a domain in $\mathbf{C}.$ Then there is a unitary
polynomial $p\in\mathcal{O}(U)[z]$ such that $E=\{(x,z)\in
U\times\mathbf{C}:p(x,z)=0\}$ and the discriminant $\Delta_p$ of
$p$ is not identically zero. $p$ will be called the optimal
polynomial for $E.$ It holds: $\tilde{\pi}(\mathcal{S}(E))=\{x\in
U:\Delta_p(x)=0\},$ where
$\tilde{\pi}:U\times\mathbf{C}\rightarrow U$ is the natural
projection.

Finally, for any analytic subset $X$ of an open set
$\tilde{U}\subset\mathbf{C}^m$ let $X_{(k)}\subset\tilde{U}$
denote the union of all irreducible components of $X$ of dimension
$k.$

\section{Approximation}\label{proofs}
\subsection{Proof of Theorem \ref{main}}\label{mainIprof}
Theorem \ref{main}, formulated in Section \ref{introduction}, is
equivalent to its slightly stronger version in which $Q$ is a Nash
mapping defined on some neighborhood of the graph of $f.$ Both
versions are in turn equivalent to the following
\begin{proposition}\label{claglo}Let $U, V$ be a domain in $\mathbf{C}^n$ and an algebraic subvariety of
$\mathbf{C}^{\hat{m}}$ respectively. Let $F:U\rightarrow V$ be a
holomorphic mapping. Then for every $x_0\in U$ there are an open
neighborhood $U_0\subset U$ and a sequence
$\{F^{\nu}:U_0\rightarrow V\}$ of Nash mappings converging
uniformly to $F|_{U_0}.$
\end{proposition}

The fact that Proposition \ref{claglo} does imply (the stronger
version of) Theorem~\ref{main} is well known (see \cite{Lem}) and
we prove it below only for completeness. The converse implication
is clear.

Let $f:U\rightarrow\mathbf{C}^k$ be the holomorphic mapping from
Theorem \ref{main}. Assume that $Q$ is a Nash mapping defined on
some neighborhood $\hat{U}$ of the graph of $f$ in
$\mathbf{C}^n\times\mathbf{C}^k,$ such that $Q(x,f(x))=0$ for
$x\in U.$ Next put $F(x)=(x,f(x))$ and $\hat{m}=n+k.$ Let $V$ be
the intersection of all algebraic subvarieties of
$\mathbf{C}^{\hat{m}}$ containing $F(U).$ Then by Proposition
\ref{claglo} there is a sequence $\{F^{\nu}:U_1\rightarrow V\}$ of
Nash mappings converging to $F|_{U_1},$ where $U_1\subset U$ is a
neighborhood of a fixed $x_0.$

We need to show that the first $n$ components of $F^{\nu}$ may be
assumed to constitute the identity and that $Q\circ F^{\nu}=0$ for
sufficiently large $\nu$. To this end denote
$Y=\{(x,v)\in\hat{U}\subset\mathbf{C}^n\times\mathbf{C}^k:Q(x,v)=0\}.$
Clearly, we may assume that $F^{\nu}(U_1)\subset\hat{U}$ for
almost all $\nu.$ Next observe that $F^{\nu}(U_1)\subset Y$ for
almost all $\nu.$ Indeed, take $\hat{z}\in F(U_1)\cap Reg(V)$ (the
intersection is non-empty as $F(U_1)\subset Sing(V)$ implies, by
the connectedness of $U,$ that $F(U)\subset Sing(V)\varsubsetneq V
$). Let $B$ be an open neighborhood of $\hat{z}$ in
$\mathbf{C}^n\times\mathbf{C}^k$ such that $B\cap V$ is a
connected manifold and let $U_2$ be a non-empty open subset of
$U_1$ such that $F(U_2), F^{\nu}(U_2)\subset B$ for almost all
$\nu.$ Then $B\cap V\subset Y$ (otherwise $F(U_2)\subset\tilde{V}$
where $\tilde{V}$ is an algebraic subvariety of
$\mathbf{C}^n\times\mathbf{C}^k$ with $dim(\tilde{V})<dim(V)$).
This implies that $F^{\nu}(U_2)\subset Y$ for almost all $\nu$
hence $F^{\nu}(U_1)\subset Y$ because $U_1$ is connected.

Let $\tilde{F}^{\nu}:U_1\rightarrow\mathbf{C}^n,$ for
$\nu\in\mathbf{N},$ be the mapping whose components are the first
$n$ components of $F^{\nu}.$ Take a neighborhood
$U_0\subset\subset U_1$ of $x_0.$ Since $\{\tilde{F}^{\nu}\}$
converges uniformly to the identity on $U_1$ and
$U_0\subset\subset U_1$ there is a sequence
$H^{\nu}:U_0\rightarrow U_{1}$ of Nash mappings such that
$\tilde{F}^{\nu}\circ H^{\nu}=id_{U_0}$ if $\nu$ is large enough.
Consequently, $F^{\nu}\circ H^{\nu}(x)=(x,f^{\nu}(x))$ for $x\in
U_0$ and $\{f^{\nu}:U_0\rightarrow\mathbf{C}^k\}$ satisfies the
assertion
of Theorem \ref{main}.\vspace*{2mm}\\
\textit{Proof of Proposition \ref{claglo}.} First observe that
since $U$ is connected, $F(U)$ is contained in one irreducible
component of $V$ so we may assume that $V$ is of pure dimension,
say $m.$

We may also assume that
$V\subset\mathbf{C}^{\hat{m}}\thickapprox\mathbf{C}^m\times\mathbf{C}^{s}$
is with proper projection onto $\mathbf{C}^m.$ Indeed, for the
generic $\mathbf{C}$-linear isomorphism
$J:\mathbf{C}^{m+s}\rightarrow\mathbf{C}^{m+s}$ the image $J(V)$
is with proper projection onto $\mathbf{C}^m.$ Thus if there
exists a sequence $H^{\nu}:U_0\rightarrow J(V)$ of Nash mappings
converging to $J\circ F|_{U_0}$ then the sequence $\{J^{-1}\circ
H^{\nu}\}$ satisfies the assertion of the proposition.

Now the problem is reduced to the case where $V$ is a hypersurface
(compare \cite{vD}, p. 394). Any $\mathbf{C}$-linear form
$L:\mathbf{C}^s\rightarrow\mathbf{C}$ determines the mapping
$\Phi_{L}:\mathbf{C}^m\times\mathbf{C}^s\rightarrow\mathbf{C}^m\times\mathbf{C}$
by the formula $\Phi_L(u,v)=(u,L(v)).$ Since $V$ is an algebraic
subset of $\mathbf{C}^m\times\mathbf{C}^s$ with proper projection
onto $\mathbf{C}^m$ then $\Phi_L(V)$ is an algebraic subset of
$\mathbf{C}^m\times\mathbf{C}$ also with proper projection onto
$\mathbf{C}^m$ for every form $L.$ Take $L$ such that the
cardinalities of the fibers of the projections of $\Phi_L(V)$ and
$V$ onto $\mathbf{C}^m$ are equal over almost every point of
$\mathbf{C}^m.$ (The generic $L$ has this property.) The set
$\Phi_L(V)$ is described by the unitary polynomial $P_L$ in one
variable (corresponding to the last coordinate of
$\mathbf{C}^m\times\mathbf{C}$) whose coefficients are polynomials
in $m$ variables and whose discriminant is non-zero (see
Section~\ref{setswithprop}).

Let $f_1,\ldots,f_m,f_{m+1},\ldots,f_{m+s}$ denote the coordinates
of the mapping $F.$ We may assume that $R_L(f_1,\ldots,f_m)\neq
0,$ where $R_L$ denotes the discriminant of $P_L.$ Indeed,
otherwise we return to the very beginning of the proof with $V$
replaced by $V\cap \{R_L=0\}$ which also contains the image of
$F.$ Since the latter variety is of pure dimension $m-1,$ the
procedure must stop.

Put $\tilde{f}=L(f_{m+1},\ldots,f_{m+s}).$ Now, $F,$ $L,P_L,V$
satisfy the following lemma which will be useful to us.
\begin{lemma}\label{reducone}Let
$\{f_1^{\nu}\},\ldots,\{f_m^{\nu}\},\{\tilde{f}^{\nu}\}$ be
sequences of holomorphic functions converging locally uniformly to
$f_1,\ldots,f_m,\tilde{f}$ respectively such that
$$P_L(f_1^{\nu},\ldots,f_m^{\nu},\tilde{f}^{\nu})=0, \mbox{ for every }\nu\in\mathbf{N}.$$ Then there
exist sequences of holomorphic functions
$\{f_{m+1}^{\nu}\},\ldots,\{f_{m+s}^{\nu}\}$ converging locally
uniformly to $f_{m+1},\ldots,f_{m+s}$ respectively such that the
image of the mapping
$(f^{\nu}_1,\ldots,f^{\nu}_m,f^{\nu}_{m+1},\ldots,f^{\nu}_{m+s})$
is contained in $V.$
\end{lemma}
\textit{Proof of Lemma \ref{reducone}.} For any holomorphic
mapping $H:E\rightarrow\mathbf{C}^m,$ where $E$ is an open subset
of $\mathbf{C}^n$ and any algebraic subvariety $X$ of
$\mathbf{C}^m\times\mathbf{C}^s$ denote
$$\mathcal{V}(X,H)=\{(x,v)\in E\times\mathbf{C}^s:(H(x),v)\in X\}.$$
Next put $\Psi_L(x,v)=(x,L(v))$ for any $x\in\mathbf{C}^n,
v\in\mathbf{C}^s.$ Let
$\tilde{\pi}:\mathbf{C}^n\times\mathbf{C}\rightarrow\mathbf{C}^n,$
$\pi:\mathbf{C}^n\times\mathbf{C}^s\rightarrow\mathbf{C}^n$ denote
the natural projections. Assume the notation of Section
\ref{setswithprop}. Then the following remark is clearly true.
\begin{remark}\label{funnyobvious}\em
Let $Z\subset E \times\mathbf{C}^s$ be an analytic subset of pure
dimension $n$ with proper projection onto a domain
$E\subset\mathbf{C}^n$ such that
$s(\pi|_Z)=s(\tilde{\pi}|_{\Psi_{L}(Z)}).$ Then for every
irreducible analytic component $\Sigma$ of $\Psi_{L}(Z)$ there
exists an irreducible analytic component $\Gamma$ of $Z$ such that
$\Psi_{L}(\Gamma)=\Sigma$ and
$s(\pi|_{\Gamma})=s(\tilde{\pi}|_{\Sigma}).$
\end{remark}
The remark allows to complete the proof of Lemma \ref{reducone}.
Put $\tilde{F}=(f_1,\ldots,f_m),$
$\tilde{F}^{\nu}=(f_1^{\nu},\ldots,f_m^{\nu}),$
$G=(f_{m+1},\ldots,f_{m+s}).$ First observe that the fact that
$R_L\circ\tilde{F}\neq 0$ and the way $L$ has been chosen imply
that the cardinalities of the generic fibers in
$\Psi_{L}(\mathcal{V}(V,\tilde{F})),$ $\mathcal{V}(V,\tilde{F}),$
$\Psi_{L}(\mathcal{V}(V,\tilde{F}^{\nu}))$ and in
$\mathcal{V}(V,\tilde{F}^{\nu})$ over $U$ are equal for large
$\nu.$ Therefore we may apply Remark \ref{funnyobvious} with
$Z=\mathcal{V}(V,\tilde{F}^{\nu})$ (for large $\nu$) to obtain the
mapping $G^{\nu}:U\rightarrow\mathbf{C}^s$ such that
$graph(G^{\nu})\subset\mathcal{V}(V,\tilde{F}^{\nu})$ and $L\circ
G^{\nu}=\tilde{f}^{\nu}.$

Observe that $\{G^{\nu}\}$ converges to $G$ locally uniformly.
Indeed, since $\{G^{\nu}\}$ is locally uniformly bounded then
taking any compact subset $K$ of $U$ and passing to a subsequence
we may assume that $\{G^{\nu}\}$ has a limit $\bar{G}$ on $K$ with
$graph(\bar{G})\subset\mathcal{V}(V,\tilde{F}).$ Assumption that
$\bar{G}\neq G|_K,$ the fact that
$\Psi_L(\mathcal{V}(V,\tilde{F}))\cap(\{x\}\times\mathbf{C})$ and
$\mathcal{V}(V,\tilde{F})\cap(\{x\}\times\mathbf{C}^s)$ have the
same number of elements for almost every $x\in U$ and the facts
that $graph(\bar{G})\subset\mathcal{V}(V,\tilde{F})$ and
$L\circ\bar{G}=L\circ G$ give a contradiction.\qed\vspace*{2mm}\\
\textit{Proof of Proposition \ref{claglo} (continuation).} Without
loss of generality assume that $x_0=0\in\mathbf{C}^n.$ To complete
the proof of Proposition \ref{claglo} it is sufficient to show
that there are sequences of Nash functions
$\{f_1^{\nu}\},\ldots,\{f_m^{\nu}\},\{\tilde{f}^{\nu}\}$
converging to $f_1,\ldots,f_m,\tilde{f}$ respectively, in some
neighborhood of $0,$ such that
$P_L(f_1^{\nu},\ldots,f_m^{\nu},\tilde{f}^{\nu})=0.$ Then by Lemma
\ref{reducone} there are sequences of holomorphic functions
$\{f^{\nu}_{m+1}\},\ldots,\{f^{\nu}_{m+s}\}$ converging to
$f_{m+1},\ldots,f_{m+s}$ respectively with the image of
$(f^{\nu}_1,\ldots,f^{\nu}_m,f^{\nu}_{m+1},\ldots,f^{\nu}_{m+s})$
contained in $V.$ Therefore
$$graph(f_{m+1}^{\nu},\ldots,f_{m+s}^{\nu})\subset\mathcal{V}(V,(f_1^{\nu},\ldots,f_m^{\nu}))$$
so for every $i=1,\ldots,s$ it holds
$$graph(f^{\nu}_{m+i})\subset\pi_i(\mathcal{V}(V,(f_1^{\nu},\ldots,f_m^{\nu})))$$
where
$\pi_i:\mathbf{C}^m\times\mathbf{C}_1\times\ldots\times\mathbf{C}_s\rightarrow\mathbf{C}^m\times\mathbf{C}_i$
denotes the natural projection. Consequently, $f_{m+i}^{\nu}$ is a
Nash function as
$\pi_i(\mathcal{V}(V,(f_1^{\nu},\ldots,f_m^{\nu})))$ is a Nash
hypersurface (see Section \ref{prelnash}).

Let us show that there are sequences of Nash functions
$\{f_1^{\nu}\},\ldots,\{f_m^{\nu}\},\{\tilde{f}^{\nu}\}$
converging to $f_1,\ldots,f_m,\tilde{f}$ respectively such that
$P_L(f_1^{\nu},\ldots,f_m^{\nu},\tilde{f}^{\nu})=0$ (i.e. let us
prove Proposition \ref{claglo} in the simplified situation where
$V$ is a hypersurface and the image of $F$ is not contained in
$Sing(V)$). This will be done by induction on $n$ (the number of
the variables $F$ depends on). Below the reduction to a lower
dimensional case (by applying the Weierstrass preparation theorem)
will be carried out in a similar way as in \cite{Ar1}. The
Tougeron implicit functions theorem, which often appears in the
context, is replaced here, as in \cite{vD}, by the following
lemma.

Let $B_n(r)$ denote a compact ball in $\mathbf{C}^n$ of radius
$r.$
\begin{lemma}\label{vddlem}\em{(\cite{vD}, p. 393). }\em
Let $d$ be a positive integer and let $M,r$
be positive real numbers. There is $\varepsilon>0$ such that for
all $A=a_0z^d+a_1z^{d-1}+\ldots+a_d\in\mathcal{O}(B_n(r))[z]$ with
$\sup_{x\in B_n(r)}|a_i(x)|<M$ where $i=0,\ldots,d$ and for all
$\alpha, c\in\mathcal{O}(B_n(r))$ with $\sup_{x\in
B_n(r)}|\alpha(x)|<M,$ $\sup_{x\in B_n(r)}|c(x)|<\varepsilon$ such
that $A\circ\alpha=c\cdot(\frac{\partial A}{\partial
z}\circ\alpha)^2$ the following holds: there is
$b\in\mathcal{O}(B_n(r))$ with $A\circ b=0$ and $\sup_{x\in
B_n(r)}|b(x)-\alpha(x)|\leq 2\sup_{x\in
B_n(r)}|c(x)(\frac{\partial A}{\partial z}\circ\alpha)(x)|.$
\end{lemma}

Before the case $n=1$ is established and the induction hypothesis
is applied we need some preparations. As above put
$\tilde{F}=(f_1,\ldots,f_m).$ Since $R_L\circ\tilde{F}\neq 0$ and
$P_L(\tilde{F},\tilde{f})=0$ it holds $\frac{\partial
P_L}{\partial z}(\tilde{F},\tilde{f})\neq 0.$ By the Weierstrass
preparation theorem, after the generic linear change of the
coordinates in $\mathbf{C}^n$ there are open neighborhoods of
zeroes $E, E'$ in $\mathbf{C}^{n-1},\mathbf{C}$ respectively such
that $\frac{\partial P_L}{\partial
z}(\tilde{F},\tilde{f})(x)=\hat{H}(x)W(x),$ $x=(x',x_n)\in E\times
E'\subset\mathbf{C}^{n-1}\times\mathbf{C},$ for some non-vanishing
function $\hat{H}\in\mathcal{O}(E\times E')$ and some unitary
polynomial $W\in\mathcal{O}(E)[x_n]$ such that
$W^{-1}(0)\cap(E\times\partial E')=\emptyset.$

Dividing $f_j,\tilde{f}$ by $W$ to the power $2$ one obtains:\\
$$f_j(x)=H_j(x)W(x)^2+r_j(x),$$ $$\tilde{f}(x)=\tilde{H}(x)W(x)^2+\tilde{r}(x)$$
for $x=(x',x_n)\in E\times E',$ where
$r_j(x),\tilde{r}(x)\in\mathcal{O}(E)[x_n]$ satisfy $deg(r_j),
deg(\tilde{r})<deg(W^2)$ and $H_j,\tilde{H}\in\mathcal{O}(E\times
E')$ for $j=1,\ldots,m.$

Denote $d=deg(W).$ The polynomials $W,r_j,\tilde{r},$ for
$j=1,\ldots,m,$ are of the
form:\vspace*{2mm}\\
$W(x)=x_n^{d}+x_n^{d-1}a_{1}(x')+\ldots+a_{d}(x'),$\\
$r_j(x)=x_n^{2d-1}b_{j,0}(x')+x_n^{2d-2}b_{j,1}(x')+\ldots+b_{j,2d-1}(x'),$\\
$\tilde{r}(x)=x_n^{2d-1}c_{0}(x')+x_n^{2d-2}c_{1}(x')+\ldots+c_{2d-1}(x').$\vspace*{2mm}\\
Replacing the holomorphic coefficients
$$a_{1},\ldots,a_{d},b_{j,0},\ldots,b_{j,2d-1},
c_{0},\ldots,c_{2d-1}$$ for all $j$ in $W,r_j,\tilde{r}$ by new
variables denoted by the same letters we obtain polynomials
$C,w_j,\tilde{w}$ respectively. Define:
$$\alpha_j=C^2S_j+w_j, \beta=C^2\tilde{S}+\tilde{w}$$ for
$j=1,\ldots,m$ where $S_j,\tilde{S}$ are new variables. Now divide
$P_L(\alpha_1,\ldots,\alpha_m,\beta)$ by $C^2$ (treated as a
polynomial in $x_n$ with polynomial coefficients) and divide
$\frac{\partial P_L}{\partial z }(\alpha_1,\ldots,\alpha_m,\beta)$
by $C$ to
obtain\vspace*{2mm}\\
(*)\hspace*{7mm}$P_L(\alpha_1,\ldots,\alpha_m,\beta)=\tilde{W}C^2+x_n^{2d-1}T_1+x_n^{2d-2}T_2+\ldots+T_{2d},$\vspace*{2mm}\\
(**)\hspace*{5mm}$\frac{\partial P_L}{\partial z
}(\alpha_1,\ldots,\alpha_m,\beta)=\bar{W}C+x_n^{d-1}T_{2d+1}+x_n^{d-2}T_{2d+2}+\ldots+T_{3d},$
where $\tilde{W},\bar{W},T_1,\ldots,T_{3d}$ are polynomials such
that $T_1,\ldots,T_{3d}$ depend only on the variables
$$a_{1},\ldots,a_{d},b_{j,0},\ldots,b_{j,2d-1},
c_{0},\ldots,c_{2d-1}.$$ Let
$$a_{1},\ldots,a_{d},b_{j,0},\ldots,b_{j,2d-1},
c_{0},\ldots,c_{2d-1},$$ where $j=1,\ldots,m,$ again denote the
holomorphic coefficients of $W, r_1,\ldots,r_m, \tilde{r}$ as at
the beginning. The tuple consisting of all these functions
satisfies the system of equations $T_1=\ldots=T_{3d}=0$ and may be
uniformly approximated in some neighborhood of
$0\in\mathbf{C}^{n-1}$ by a sequence (indexed by $\nu$) of tuples
of Nash functions
$$a_{1}^{\nu},\ldots,a_{d}^{\nu},b_{j,0}^{\nu},\ldots,b_{j,2d-1}^{\nu},
c_{0}^{\nu},\ldots,c_{2d-1}^{\nu},$$ where $j=1,\ldots,m,$ also
satisfying the system $T_1=\ldots=T_{3d}=0.$ Indeed, if $n=1$ then
the coefficients are constant and therefore may be taken as their
own approximations. If $n>1$ we are done by the induction
hypothesis because the coefficients depend on $n-1$ variables.

Using the obtained Nash functions define:\vspace*{2mm}\\
$W_{\nu}(x)=x_n^{d}+x_n^{d-1}a_{1}^{\nu}(x')+\ldots+a_{d}^{\nu}(x'),$\\
$r_{j,\nu}(x)=x_n^{2d-1}b_{j,0}^{\nu}(x')+x_n^{2d-2}b_{j,1}^{\nu}(x')+\ldots+b_{j,2d-1}^{\nu}(x'),$\\
$\tilde{r}_{\nu}(x)=x_n^{2d-1}c_{0}^{\nu}(x')+x_n^{2d-2}c_{1}^{\nu}(x')+\ldots+c_{2d-1}^{\nu}(x').$\vspace*{2mm}\\
Finally for, $j=1,\ldots,m,$ define
$$f_j^{\nu}(x)=H_j^{\nu}(x)(W_{\nu}(x))^2+r_{j,\nu}(x),$$ $$\bar{f}^{\nu}(x)=\tilde{H}^{\nu}(x)(W_{\nu}(x))^2+\tilde{r}_{\nu}(x).$$
Here $H^{\nu}_{j},\tilde{H}^{\nu},$ are any sequences of
polynomials converging uniformly to $H_{j},\tilde{H}$ in some
neighborhood of $0\in\mathbf{C}^n$ respectively.

Now it is easy to see that by (*), (**) and the way
$f^{\nu}_1,\ldots,f^{\nu}_m,\bar{f}^{\nu}$ are defined, there is a
neighborhood of $0\in\mathbf{C}^n$ in which for every $\nu$ the
following holds:
$$P_L(f_1^{\nu},\ldots,f_m^{\nu},\bar{f}^{\nu})=R^{\nu}(\frac{\partial{P_L}}{\partial{z}}(f_1^{\nu},\ldots,f_m^{\nu},\bar{f}^{\nu}))^2,$$
where $\{R^{\nu}\}$ is a sequence of holomorphic functions
converging to zero. Therefore it suffices to apply Lemma
\ref{vddlem} with $A=P_L(f_1^{\nu},\ldots,f_m^{\nu},z),$
$\alpha=\bar{f}^{\nu}$ and $c=R^{\nu}$ (for sufficiently large
$\nu$) to obtain
$$P_L(f_1^{\nu},\ldots,f_m^{\nu},\tilde{f}^{\nu})=0$$ in some neighborhood $U_0$ of zero for
every $\nu\in\mathbf{N},$ where $\{\tilde{f}^{\nu}\}$ is a
sequence of holomorphic functions converging to $\tilde{f}$ in
$U_0.$\qed
%
\subsection{Algorithm}\label{examp}Section \ref{examp} is devoted
to a recursive algorithm of Nash appro\-ximation of a holomorphic
mapping $F:U\rightarrow V\subset\mathbf{C}^{\hat{m}},$ where $U$
is a neighborhood of zero in $\mathbf{C}^n$ and $V$ is an
algebraic variety. The correctness of the algorithm follows from
the proof of Proposition \ref{claglo}.

The method of approximation presented here is more efficient than
the one we developed in \cite{B6}. The main difference is that now
the polynomial $R_L$ defined in step 4 below need not be
factorized into powers of pairwise distinct optimal polynomials.

For $\nu\in\mathbf{N},$ the approximating mapping
$F^{\nu}=(f^{\nu}_1,\ldots,f^{\nu}_{\hat{m}}):U_0\rightarrow V,$
returned as the output of the algorithm, is represented by
$\hat{m}$ non-zero polyno\-mials
$P_i^{\nu}(x,z_i)\in(\mathbf{C}[x])[z_i],$ $i=1,\ldots,\hat{m},$
such that $P_i^{\nu}(x,f^{\nu}_i(x))=0$ for $x\in U_0.$ We
restrict attention to the local case i.e. $U_0$ is an open
neighborhood of a fixed $x_0\in U.$
More precisely, we work with the following data:\vspace*{2mm}\\
\textbf{Input:} a holomorphic mapping
$F=(f_1,\ldots,f_{\hat{m}}):U\rightarrow
V\subset\mathbf{C}^{\hat{m}},$ $F=F(x),$ where $U$ is an open
neighborhood of $0\in\mathbf{C}^n$ and $V$ is an algebraic
variety.\vspace*{2mm}\\
\textbf{Output:} $P_i^{\nu}(x,z_i)\in(\mathbf{C}[x])[z_i],$
$P_i^{\nu}\neq 0$ for $i=1,\ldots,\hat{m}$ and $\nu\in\mathbf{N},$
with the following properties:\vspace*{1mm}\\
(a) $P_i^{\nu}(x,f_i^{\nu}(x))=0$ for every $x\in U_0,$ where
$F^{\nu}=(f^{\nu}_1,\ldots,f^{\nu}_{\hat{m}}):U_0\rightarrow V$ is
a holomorphic mapping such that $\{F^{\nu}\}$ converges uniformly
to $F$ on an open neighborhood $U_0$ of
$0\in\mathbf{C}^n,$\\
(b) $P_i^{\nu}$ is a unitary polynomial in $z_i$ of degree
independent of $\nu$ whose co\-efficients (belonging to
$\mathbf{C}[x]$) converge uniformly to holomorphic functions on
$U_0$ as $\nu$ tends to infinity.\vspace*{2mm}

First let us comment on the notation and the idea of the
algorithm. The meaning of the symbol $V_{(m)}$ and the notion of
the optimal polynomial used below can be found in Subsection
\ref{setswithprop}.

In steps 2 and 6 we apply linear changes of the coordinates.
Having approximated the mapping $\hat{J}\circ F\circ
J|_{J^{-1}(U)}:J^{-1}(U)\rightarrow\hat{J}(V),$ where
$\hat{J}:\mathbf{C}^{\hat{m}}\rightarrow\mathbf{C}^{\hat{m}},$
$J:\mathbf{C}^n\rightarrow\mathbf{C}^n$ are linear isomorphisms,
one can obtain the output data for $F$ following standard
arguments. (Composing $F$ and $J$ does not lead to any
difficulties. As for $\hat{J},$ it is sufficient to use the fact
that the integral closure of a commutative ring in another
commutative ring is again a ring.) Therefore, when the coordinates
are changed, we write what (as a result) may be assumed about the
mapping $F,$ but the notation is left unchanged.

The aim of steps 1-3 is to prepare the variety $V$ so that for the
polynomial $P_L$ calculated in step 4, Lemma \ref{reducone} holds
(for details see the proof of Proposition \ref{claglo}). Steps 5-9
are responsible for the fact that
$P_L(f_1^{\nu},\ldots,f_m^{\nu},\bar{f}^{\nu})=R^{\nu}(\frac{\partial{P_L}}{\partial{z}}(f_1^{\nu},\ldots,f_m^{\nu},\bar{f}^{\nu}))^2,$
where $f_1^{\nu},\ldots,f_m^{\nu}$ are defined in step 10 whereas
$\{R^{\nu}\}$ is a sequence of holomorphic functions converging to
zero and $\{\bar{f}^{\nu}\}$ is a sequence of holomorphic
functions converging to $L(f_{m+1},\ldots,f_{m+s}).$ This enables
to use Lemma \ref{vddlem} which together with Lemma \ref{reducone}
implies that $P^{\nu}_{m+1},\ldots,P^{\nu}_{m+s}$ are well defined
in step 11. As for $P^{\nu}_{1},\ldots,P^{\nu}_{m},$ these
polynomials are obtained in step 10 by applying the results of the
algorithm switched for
the lower dimensional case in step 9.\vspace*{2mm}\\
\textbf{Algorithm:} \textbf{1.} Replace $V$ by $V_{(m)}$ such that
$F(U)\subset V_{(m)}.$\\
\textbf{2.} Apply a linear change of the coordinates in
$\mathbf{C}^{\hat{m}}$ after which $\rho|_{V}$ is a proper
mapping, where
$\rho:\mathbf{C}^m\times\mathbf{C}^s\approx\mathbf{C}^{\hat{m}}\rightarrow\mathbf{C}^m$
is the natural projection.\\
\textbf{3.} Choose a $\mathbf{C}$-linear form
$L:\mathbf{C}^s\rightarrow\mathbf{C}$ such that the generic fibers
of $\rho|_{V}$ and $\tilde{\rho}|_{\Phi_L(V)}$ over $\mathbf{C}^m$
have the same cardinalities. Here
$\tilde{\rho}:\mathbf{C}^m\times\mathbf{C}\rightarrow\mathbf{C}^m$
is the natural projection and $\Phi_L(y,v)=(y,L(v))$ for
$(y,v)\in\mathbf{C}^m\times\mathbf{C}^s.$\\
\textbf{4.} Calculate the optimal polynomial
$P_L(y,z)\in(\mathbf{C}[y])[z]$ describing
$\Phi_L(V)\subset\mathbf{C}^m_y\times\mathbf{C}_z.$ Calculate the discriminant $R_L\in\mathbf{C}[y]$ of $P_L.$\\
\textbf{5.} If $R_L(f_{1},\ldots,f_{m})=0$ then return to step 2
with $m,$ $s,$ $V$ replaced by $m-1,$ $s+1,$ $V\cap\{R_L=0\}$
respectively. Otherwise put $\tilde{f}=L(f_{m+1},\ldots,f_{m+s})$
and observe that $\frac{\partial P_L}{\partial z
}(f_1,\ldots,f_m,\tilde{f})\neq 0.$\\
\textbf{6.} Apply a linear change of the coordinates in
$\mathbf{C}^n$ after which the following holds: $\frac{\partial
P_L}{\partial z }(f_1,\ldots,f_m,\tilde{f})(x)$$=\hat{H}(x)W(x)$
in some neighborhood of $0\in\mathbf{C}^n,$ where $\hat{H}$ is a
holomorphic function, $\hat{H}(0)\neq 0$ and $W$ is a unitary
polynomial in $x_n$ with holomorphic coefficients depending on
$x'=(x_1,\ldots,x_{n-1}).$ Put $d=deg(W).$\\
\textbf{7.} Divide $f_i(x),\tilde{f}(x)$ by $(W(x))^2$ to obtain
$f_i(x)=(W(x))^2H_i(x)+r_i(x)$ and
$\tilde{f}(x)=(W(x))^2\tilde{H}(x)+\tilde{r}(x)$ in some
neighborhood of $0\in\mathbf{C}^n,$ $i=1,\ldots,m.$ Here
$H_i,\tilde{H}$ are holomorphic functions and $r_i,\tilde{r}$ are
polynomials in
$x_n$ with holomorphic coefficients depending on $x',$ such that $deg(r_i), deg(\tilde{r})<2d.$\\
\textbf{8.} Treating $H_i,\tilde{H},$ $i=1,\ldots,m,$ and all the
coefficients of $W,r_1,\ldots,r_m,\tilde{r}$ as new variables
(except for the coefficient $1$ standing at the leading term of
$W$) apply the division procedure for
polynomials to obtain:\\
$P_L(W^2H_1+r_1,\ldots,W^2H_m+r_m,W^2\tilde{H}+\tilde{r})=$\\
$=\tilde{W}W^2+x_n^{2d-1}T_1+x_n^{2d-2}T_2+\ldots+T_{2d},$\\
$\frac{\partial P_L}{\partial
z}(W^2H_1+r_1,\ldots,W^2H_m+r_m,W^2\tilde{H}+\tilde{r})=$\\
$=\bar{W}W+x_n^{d-1}T_{2d+1}+x_n^{d-2}T_{2d+2}+\ldots+T_{3d}.$\\
Here $T_1,\ldots,T_{3d}$ are polynomials depending only on the
variables standing for the coefficients of $W,$
$r_1,\ldots,r_m,\tilde{r}.$ Moreover, $T_1\circ
g=\ldots=T_{3d}\circ g=0,$ where $g$ is the mapping whose
components are all these coefficients
(cp. Section \ref{mainIprof}).\\
\textbf{9.} If $n>1$ then apply the Algorithm with $F, V$ replaced
by $g,$$\{T_1=\ldots=T_{3d}=0\}$ respectively. As a result, for
every component $c(x')$ of $g(x')$ and for every
$\nu\in\mathbf{N}$ one obtains a unitary polynomial
$Q_c^{\nu}(x',t_c)\in(\mathbf{C}[x'])[t_c]$ which put in place of
$P_i^{\nu}(x,z_i)$ satisfies (a) and (b) above with $x,$ $z_i,$
$f_i^{\nu}$ replaced by $x',$ $t_c,$ $c^{\nu}$ respectively. Here,
for every $c,$ $\{c^{\nu}\}$ is a sequence of Nash functions
converging to $c,$ in some neighborhood of $0\in\mathbf{C}^{n-1},$
such that for every fixed $\nu$ the mapping $g^{\nu}$ obtained by
replacing every $c$ of $g$ by $c^{\nu}$
satisfies $T_1\circ g^{\nu}=\ldots=T_{3d}\circ g^{\nu}=0.$\\
If $n=1$ then $g$ is
constant and then it is its own approximation yielding the $Q^{\nu}_c$'s immediately.\\
\textbf{10.} Approximate $H_i$ for $i=1,\ldots,m,$ by a sequence
$\{H_i^{\nu}\},$ of polynomials. Let $W_{\nu},r_{1,{\nu}},$
$\ldots,r_{m,{\nu}},\tilde{r}_{\nu},$ for every
$\nu\in\mathbf{N},$ be the polynomials in $x_n$ defined by
replacing the coefficients of $W,$ $r_1,\ldots,r_m,\tilde{r}$ by
their Nash approxima\-tions (i.e. the components of $g^{\nu}$)
determined in step 9. Using $Q_c^{\nu}$ (for all $c$) and
$H^{\nu}_i$ one can calculate $P_i^{\nu}\in(\mathbf{C}[x])[z_i],$
for $i=1,\ldots,m,$ satisfying (b) and (a) with
$f_i^{\nu}=H_i^{\nu}(W_{\nu})^2+r_{i,{\nu}}$ being the $i$'th
component of the mapping $F^{\nu}$ (whose last $\hat{m}-m$
components are determined by
$P_{m+1}^{\nu},\ldots,P_{\hat{m}}^{\nu}$ obtained in the next
step). To calculate $P^{\nu}_1,\ldots,P^{\nu}_m$ one can follow
the standard proof of the fact that the integral closure of a
commutative ring in another commutative ring is again a ring.\\
\textbf{11.} Put
$V^{\nu}=\{(x,z)\in\mathbf{C}^n_x\times\mathbf{C}^{m+s}_{z} : z\in
V, P^{\nu}_i(x,z_i)=0 \mbox{ for } i=1,\ldots,m\},$ where
$z=(z_1,\ldots,z_m,z_{m+1},\ldots,z_{m+s}).$ For $i=1,\ldots,s$
and $\nu\in\mathbf{N}$ take
$P^{\nu}_{m+i}\in(\mathbf{C}[x])[z_{m+i}]$ to be the optimal
polynomial describing the image of the projection of $V^{\nu}$
onto $\mathbf{C}^n_x\times\mathbf{C}_{z_{m+i}}.$


\begin{thebibliography}{}
\bibitem{An}
Andr\'e, M.: Cinq expos\'es sur la d\'esingularization,
manuscript, \'Ecole Polytechnique F\'ed\'erale de Lausanne, 1992
\bibitem{Ar1}
Artin, M.: On the solutions of analytic equations. Invent. Math.
\textbf{5}, 277-291 (1968)
\bibitem{Ar2}
Artin, M.: Algebraic approximation of structures over complete
local rings. Publ. I.H.E.S. \textbf{36}, 23-58 (1969)
\bibitem{Ar3}
Artin, M.: Algebraic structure of power series rings. Contemp.
Math. \textbf{13}, 223-227 (1982)
\bibitem{B1}
Bilski, M.: On approximation of analytic sets. Manuscripta Math.
\textbf{114}, 45-60 (2004)
\bibitem{B2}
Bilski, M.: Approximation of analytic sets with proper projection
by Nash sets. C.R. Acad. Sci. Paris, Ser. I \textbf{341}, 747-750
(2005)
\bibitem{B3}
Bilski, M.: Approximation of analytic sets by Nash tangents of
higher order. Math. Z., \textbf{256}, 705-716 (2007)
\bibitem{B6}
Bilski, M.: Algebraic approximation of analytic sets and mappings.
Preprint server: http://arxiv.org
\bibitem{Ch}
Chirka, E. M.: Complex analytic sets. Kluwer Academic Publ.,
Dordrecht-Boston-London 1989
\bibitem{CRS}
Coste, M., Ruiz, J. M., Shiota, M.: Approximation in compact Nash
manifolds. Amer. J. Math. \textbf{117}, 905-927 (1995)
\bibitem{vD}
van den Dries, L.: A specialization theorem for analytic functions
on compact sets. Nederl. Akad. Wetensch. Indag. Math. \textbf{44},
391-396 (1982)
\bibitem{Lem}
Lempert, L.: Algebraic approximations in analytic geometry.
Invent. Math. \textbf{121}, 335-354 (1995)
\bibitem{Og}
Ogoma, T.: General N\'eron desingularization based on the idea of
Popescu. J. of Algebra \textbf{167}, 57-84 (1994)
\bibitem{Po1}
Popescu, D.: General N\'eron desingularization. Nagoya Math. J.
\textbf{100}, 97-126 (1985)
\bibitem{Po2}
Popescu, D.: General N\'eron desingularization and approximation.
Nagoya Math. J. \textbf{104}, 85-115 (1986)
\bibitem{Sp}
Spivakovsky, M.: A new proof of D. Popescu's theorem on smoothing
of ring homomorphisms. J. Amer. Math. Soc., \textbf{12}, 381-444
(1999)
\bibitem{Tw}
Tworzewski, P.: Intersections of analytic sets with linear
subspaces. Ann. Sc. Norm. Super. Pisa \textbf{17,} 227-271 (1990)
\bibitem{Wh}
Whitney, H.: Complex Analytic Varieties. Addison-Wesley Publishing
Co., Reading, Mass.-London-Don Mills, Ont., 1972
\end{thebibliography}
\end{document}